\documentclass[11pt]{amsart}

\usepackage[english]{babel}
\usepackage[utf8]{inputenc}
\usepackage{amsmath}
\usepackage{amssymb}
\usepackage{amsfonts}
\usepackage{amsthm}
\usepackage{mathrsfs}
\usepackage[all]{xy}
\usepackage[pdftex]{graphicx}
\usepackage{color}
\usepackage{cite}
\usepackage{url}
\usepackage{indent first}
\usepackage[labelfont=bf,labelsep=period,justification=raggedright]{caption}
\usepackage[english]{babel}
\usepackage[utf8]{inputenc}
\usepackage{hyperref}
\usepackage[capitalize,nameinlink]{cleveref}
\usepackage[colorinlistoftodos]{todonotes}
\usepackage{tkz-fct}
\usepackage[margin=1in]{geometry} 
\usepackage{tikz}
\usetikzlibrary{decorations.markings}

\newcommand{\eps}{\epsilon}

\newcommand{\EE}{\mathbb{E}}

\newcommand{\NN}{\mathbb{N}}

\newcommand{\mcG}{\mathcal{G}}

\newcommand{\ind}{\mathbf{1}}

\theoremstyle{plain}
\newtheorem{thm}{Theorem}
\newtheorem{lemma}[thm]{Lemma}
\newtheorem{cor}[thm]{Corollary}
\newtheorem{conj}[thm]{Conjecture}
\newtheorem{prop}[thm]{Proposition}

\newtheorem{qn}[thm]{Question}
\newtheorem{claim}[thm]{Claim}
\newtheorem{obs}[thm]{Observation}

\theoremstyle{definition}
\newtheorem{defn}[thm]{Definition}

\theoremstyle{remark}
\newtheorem{rem}[thm]{Remark}

\numberwithin{equation}{section}
\numberwithin{thm}{section}

\begin{document}

\title{Variations on Sidorenko's conjecture in tournaments}
\author[Fox]{Jacob Fox}
\author[Himwich]{Zoe Himwich}
\author[Mani]{Nitya Mani}
\author[Zhou]{Yunkun Zhou}

\thanks{Fox was supported by NSF Award DMS-2154129.}
\thanks{Himwich was supported by the Fernholz Foundation Minerva Fellowship Program.}
\thanks{Mani was supported by the NSF Graduate Research Fellowship Program and a Hertz Graduate Fellowship.}
\thanks{Zhou was supported by NSF GRFP Grant DGE-1656518.}

\address{Fox \& Zhou: Department of Mathematics, Stanford University, Stanford, CA, USA}
\email{\{jacobfox,yunkunzhou\}@stanford.edu}

\address{Himwich: Department of Mathematics, Columbia University, New York, NY, USA}
\email{zmchimwich@gmail.com}

\address{Mani: Department of Mathematics, Massachusetts Institute of Technology, Cambridge, MA 02139, USA}
\email{nmani@mit.edu}

\date{\today}

\begin{abstract} 
We study variants of Sidorenko's conjecture in tournaments, where new phenomena arise that do not have clear analogues in the setting of undirected graphs. We first consider oriented graphs that are systematically under-represented in tournaments (called \textit{tournament anti-Sidorenko}). We prove that such oriented graphs must be quite sparse; specifically, the maximum number of edges of a $k$-vertex oriented graph which is tournament anti-Sidorenko is $(1+o(1))k\log_2 k$. 

We also give several novel constructions of oriented graphs that are systematically over-represented in tournaments (\textit{tournament Sidorenko}); as a representative example, we show that most ways to delete an edge from a transitive tournament yield a tournament Sidorenko oriented graph. As an illustration of our methods, we characterize which  orientations of stars are tournament Sidorenko and which are tournament anti-Sidorenko. 
\end{abstract}

\maketitle

\section{Introduction}
\subsection{Background}
A fundamental problem in extremal graph theory is estimating the number of copies of a graph $H$ in a graph $G$ with a given number of vertices and edges.
The \textit{$H$-density} of a graph $G = (V(G), E(G))$ on $v(G)$ vertices and $e(G)$ edges is the fraction of injective vertex maps $\rho: V(H) \rightarrow V(G)$ that send edges to edges. Each $\rho$ gives a distinct labeled copy of $H$ in $G$. It is possible to obtain an upper bound on the minimum $H$-density over all graphs $G$ on $n$ vertices by computing the expected $H$ density, when $G = \mcG(n, p)$ is an Erd\H{o}s-Renyi random graph. 

A motivating problem in this area of graph theory is \textit{Sidorenko's conjecture}, the proposition that the above upper bound on the minimal $H$-density is sharp when $H$ is bipartite. 
This conjecture, independently posed by Erd\H{o}s-Simonovits \cite{SIM82} and Sidorenko \cite{SID93},
is often stated in the language of \textit{graph homomorphisms}, vertex maps of graphs that send edges to edges.

\begin{defn}\label{d:homdensity}
The \textit{homomorphism density} of $H$ in $G$ is the fraction of vertex maps $H \rightarrow G$ that are homomorphisms, given by $t_H(G) = \frac{h_H(G)}{v(G)^{v(H)}},$ where $h_H(G)$ is the number of homomorphisms $H \rightarrow G$.
\end{defn}

\begin{conj}[Sidorenko's Conjecture]\label{c:sidconj}
For every bipartite undirected graph $H$ and every undirected graph $G$,
$$t_H(G) \ge t_{K_2} (G)^{e(H)}.$$
\end{conj}

Sidorenko showed~\cite{SID93} that this conjecture holds for several types of graphs, including complete bipartite graphs and even cycles. The general conjecture remains open, but a variety of special cases have been resolved. For example, Conlon, Fox, and Sudakov~\cite{CON10} showed that the conjecture holds for bipartite graphs with one vertex complete to the other side. In several works, including ~\cite{LI17,KIM16,CON18,SZE}, further bipartite subfamilies have been shown to satisfy the conjecture. A graph which satisfies Sidorenko's conjecture has the \textit{Sidorenko property}.

The Sidorenko property is closely tied to a stronger graph-theoretic property, characterizing   \textit{quasirandomness} of graphs by subgraph counts. Quasirandom graphs were first studied by Thomason \cite{TH87} and Chung-Graham-Wilson \cite{FAN89}. Both of these works observed that a large number of properties satisfied by Erd\H{o}s-Renyi random graphs are actually equivalent. 
Quasirandomness is a property defined for a family of graphs on $n$ vertices which indicates that it satisfies, with high probability as $n\to \infty$, several equivalent properties typically associated with an Erd\H{o}s-Renyi random graph.

\begin{defn}
A sequence ($G_n : n = 1, 2, \ldots$) of undirected graphs is called \textit{quasirandom with density $p$} (where $0 < p < 1$) if, for
every graph $H$, $$t_H (G_n) = (1 + o(1))p^{e(H)}.$$
\end{defn}

In the quasirandom setting, there is a stronger analogue of Sidorenko's conjecture. A graph $H$ is \textit{forcing} if a family of graphs $\{G_n\}_{n = 1}^{\infty}$ is quasirandom with density $p$ if and only if the number of copies of $H$ in $G_n$ is asymptotically the number achieved in Erd\H{o}s-Renyi graphs of density $p$. The \textit{forcing conjecture}, initially posed by Skokan and Thoma~\cite{SKO04} states that subgraphs $A$ are forcing if and only if they are bipartite and contain a cycle (it is straightforward to show that these conditions are necessary).

 While there has been an extensive effort over the past decades to resolve parts of Sidorenko's conjecture and to study related extremal properties of graphs, the analogues of these problems for oriented graphs are less understood, even when restricted to the setting of \textit{tournaments,} orientations of complete graphs (i.e. oriented graphs that contain a single directed edge between every pair of distinct vertices). 

 \begin{rem}
  Throughout this article, we consider all directed graphs (abbreviated \textit{digraphs}), to be \textit{oriented}, i.e. orientations of undirected graphs without self-loops or antiparallel edges.
  \end{rem}

\subsection{Tournament anti-Sidorenko}
The story of Sidorenko-style properties in tournaments is much more nuanced than it might initially appear. For example, let $T$ be an arbitrary tournament on $n$ vertices.
There exist digraphs $D$ such that the number of copies of $D$ in $T$ is always at most the expected number of copies of $D$ in an $n$-vertex tournament with edges oriented uniformly at random. For example, if $D = \vec P_k$, the directed path with $k$ edges (where every vertex has both in- and out-degree at most $1$), then the number of copies of $D$ in $T$ is always at most $n(n/2)^k$ as observed in~\cite{SSZ20}. We thus consider the following \textit{two} Sidorenko-style properties in the context of tournaments (where homomorphism density is defined in the same fashion as it is for undirected graphs).

\begin{defn}[Tournament Sidorenko properties]\label{d:tsidorenko}
A digraph $D$ satisfies the \textit{tournament Sidorenko} property if, for all tournaments $T$, 
$$t_D(T) \ge (1-o(1))2^{-e(D)}.$$
Likewise, $D$ satisfies the \textit{tournament anti-Sidorenko} property if, for all tournaments $T$, 
$$t_D(T) \le 2^{-e(D)}.$$
\end{defn}
\begin{rem}
The asymmetry in the definitions of the tournament Sidorenko and tournament anti-Sidorenko properties is purely a consequence of our restriction to tournaments without self-loops. If we allow our tournaments to have self-loops (or work with the tournament analogue of \textit{graphons}), the $o(1)$ in the definition of the tournament Sidorenko property will disappear.
\end{rem}

In some literature (e.g.~\cite{SSZ20}), subgraphs with the tournament anti-Sidorenko property are also referred to as \textit{negative}. In addition to directed paths, directed cycles of length $r$ have the tournament anti-Sidorenko property if and only if $r \not \equiv 0 \pmod 4$, as shown in~\cite{GKLV20} as does the following orientation of an undirected $4$-cycle: $\{(a, b), (b, c), (c, d), (a, d)\}$~\cite{GRI13}.

One of the primary contributions of this article is a better understanding of which subgraphs can have the tournament anti-Sidorenko property. \cref{c:notanti} gives an asymptotically tight bound on the maximum number of edges in a tournament anti-Sidorenko digraph. 

\begin{thm}  \label{c:notanti}
Let $f(k)$ be the maximum number of edges in a $k$-vertex tournament anti-Sidorenko digraph. Then,
$$k \log_2 k - O(k) \le f(k) \le k \log_2 k.$$
\end{thm}

In \cref{s:antisid}, we also give multiple constructions of families of anti-Sidorenko digraphs which achieve this lower bound. 

Often, we can obtain examples of tournament anti-Sidorenko digraphs by building on smaller digraphs known to be systematically over or under-represented. For example, following an argument in~\cite{ZHAO19} and observation of~\cite{SSZ20}, we can obtain a new digraph with the tournament anti-Sidorenko property by taking a disjoint union of two identical copies of some digraph with the tournament anti-Sidorenko property and adding a single edge connecting a matching pair of vertices. In this vein, we give several other recursive constructions of novel families of digraphs with the tournament anti-Sidorenko property (c.f.~\cref{s:antisid}). 

\subsection{Tournament Sidorenko} We also consider digraphs with the tournament Sidorenko property. In this context, \textit{transitivity} plays a role similar to that played by the bipartite property in the undirected setting.

\begin{defn}
A digraph $D
= (V, E)$ is \textit{transitive} if, for all pairs of edges $(x, y), (y, z) \in E$, if there is an edge on $\{x, z\}$ in $E$, then it is oriented so that $(x, z) \in E$.
\end{defn}

Transitivity is a necessary condition for a digraph to have the tournament Sidorenko property, as a transitive tournament has no copies of any oriented subgraph that is not transitive.

Coregliano and Razborov~\cite{COR17} showed that transitive tournaments have the tournament Sidorenko property. It is a consequence of a more general reduction in~\cite{genFHMZ22} that any digraph with a homomorphism to an edge has the tournament Sidorenko property if its underlying undirected graph has the asymmetric Sidorenko property.

In this work, we expand the family of digraphs known to have the tournament Sidorenko property by illustrating properties that allow us to bootstrap from the known transitive digraphs that have the tournament Sidorenko property to bigger families of digraphs. A special case of a recursive construction (given in full generality in~\cref{p:growtsid}) is the following.

\begin{prop} \label{c:growspecial}
If $D = (V, F)$ has the tournament Sidorenko property, then the digraph $D' = (V \cup \{v\}, F')$ where $F' = F \cup \{(v, w) : w \in V\}$ has the tournament Sidorenko property.
\end{prop}

\begin{figure}[t]
\centering
\begin{tikzpicture}
\draw[fill=black] (0,0) circle (3pt);
\draw[fill=black] (4,0) circle (3pt);
\draw[fill=black] (2,-1) circle (3pt);
\draw[fill=blue] (2,2) circle (3pt);
\node (A) at (-0.5,0) {$a$};
\node (C) at (4.5,0) {$c$};
\node (B) at (2,-1.5) {$b$};
\node (D) at (2,2.5) {$v$};
\begin{scope}[very thick,decoration={
    markings,
    mark=at position 0.5 with {\arrow{>}}}
    ] 
    \draw[postaction={decorate}] (0,0)--(2,-1);
    \draw[postaction={decorate}] (4,0)--(2,-1);
    \draw[blue,postaction={decorate}] (2,2)--(0,0);
    \draw[blue,postaction={decorate}] (2,2)--(2,-1);
    \draw[blue,postaction={decorate}] (2,2)--(4,0);
\end{scope}
\end{tikzpicture}
\caption{The above four vertex graph has the tournament Sidorenko property, since the induced digraph on $\{a, b, c\}$ has the directed Sidorenko property and $v$ is complete to $\{a,b,c\}$. This is a special case of~\cref{c:growspecial}.}
\label{f:newtsid}
\end{figure}
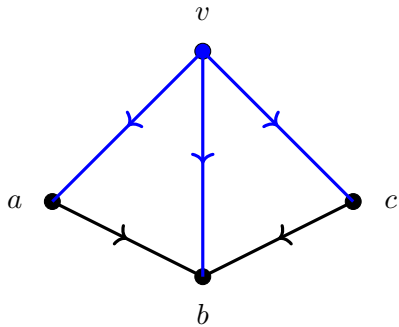

In other words, given a tournament Sidorenko digraph $D$, adding a vertex complete to $D$ and orienting each new edge the same way yields a new tournament Sidorenko digraph $D'.$ A simple example appears in \cref{f:newtsid}.
As a more complicated example,~\cref{c:growspecial} also establishes the tournament Sidorenko property for the following more interesting family of digraphs.
\begin{thm} \label{c:transminusedge}
If $T_k = ([k], E)$ is a transitive tournament with $E = \{(i, j) \in [k]^2 \mid i < j\}$ and $e = (i, j) \in E$ is any edge with $j - i \neq 2$, then $T_k \backslash \{e\}$ is tournament Sidorenko.
\end{thm}
We note that there are also ways to delete an edge from a transitive tournament that do not yield a tournament Sidorenko digraph. As an example, we can obtain a directed path of length $2$ by deleting a single edge from a transitively oriented triangle.

\subsection{Trees} While our work falls short of a complete characterization of the set of tournament Sidorenko and anti-Sidorenko digraphs, we are able to achieve a complete understanding of orientations of stars.

\begin{thm}\label{t:star}
Suppose $S$ is an orientation of an undirected star (i.e. an orientation of $K_{1, s}$).
\begin{itemize}
    \item $S$ is tournament Sidorenko if and only if it has a homomorphism to an edge
    \item $S$ is tournament anti-Sidorenko if and only if the unique vertex in $S$ of total degree $s$ has almost equal in-degree and out-degree (i.e. differing by at most one)
\end{itemize}
\end{thm}

Since undirected trees are known to have the asymmetric Sidorenko property, we more generally observe that for any undirected tree, any orientation that has a homomorphism to an edge has the tournament Sidorenko property. It is natural to wonder if every undirected tree also has an orientation that is tournament anti-Sidorenko. We conjecture that this is the case, and affirmatively resolve this conjecture for a special family of trees.

\begin{conj}\label{c:treeanti}
For every undirected tree, there exists some orientation that has the tournament anti-Sidorenko property.
\end{conj}

\begin{prop}\label{p:partialtree}
\cref{c:treeanti} holds for all undirected trees with at most one vertex of even degree.
\end{prop}

\subsection{Forcing} 
Lastly, we study an analogue of the forcing property in tournaments. A notion of quasirandommness for tournaments was introduced more than twenty years ago in~\cite{FAN91}. However, developing a robust understanding of \textit{forcing} in tournaments has been challenging, with recent exciting progress in~\cite{BU19,COR15,COR17,HA19}. We begin by giving a more precise notion of forcing for tournaments, defining \textit{quasirandom direction}.
\begin{defn}
A tournament $T$ with $n$ vertices has \textit{$\eps$-quasirandom direction} if, for all $A, B \subseteq V(T)$, we have $e(A, B) - e(B, A) \le \eps n^2$.
A sequence of tournaments $\{T_n\}_{n \in \NN}$ (where $T_n$ has $n$ vertices) has \textit{quasirandom direction} 
if each $T_n$ has $\eps_n$-quasirandom direction for some sequence $\eps_n \to 0$
 as $n \to \infty$.
 \end{defn}

 A counting lemma~\cite{FAN91}  implies that for tournament sequences $\{T_n\}$ with $\eps$-quasirandom direction, the number of copies of any fixed digraph is asymptotically what we expect of a tournament whose edges are oriented uniformly at random. For digraphs $D$ with the tournament forcing property, we have the converse, that the count of copies of $D$ in a sequence of tournaments $\{T_n\}$ is sufficient to establish that the sequence of tournaments has quasirandom direction.
\begin{defn}
Digraph $D$ has the \textit{tournament forcing property} if any sequence of tournaments $\{T_n\}_{n \in \NN}$ (where $T_n$ has $n$ vertices) has quasirandom direction if and only if for all $n \in \NN$,
$$t_D(T_n) = 1/2^{e(D)} + o(1).$$
\end{defn}

In 2019, Bucic, Long, Shapira, and Sudakov~\cite{BU19} showed that if a sequence of tournaments $\{T_n\}$ has the ``correct'' number of copies of transitive tournaments with at least $7$ vertices, any single digraph $D$ whose count certifies that $\{T_n\}$ has quasirandom direction must be transitive. Hancock et. al.~\cite{HA19} showed that apart from a single $5$ vertex non-transitive tournament, all forcing tournaments are transitive ($T_5^8$ in Figure 2 of~\cite{COR15} illustrates this non-transitive forcing tournament on $5$ vertices).

We observe that tournament forcing digraphs must fall into one of two types.
\begin{thm} \label{t:overunder}
If a digraph $D$ has the tournament forcing property, then either $D$ has the tournament Sidorenko property or $D$ has the tournament anti-Sidorenko property.
\end{thm}

This result enables us to observe that under a weak density condition, transitivity is necessary for a digraph to have the tournament forcing property.

\begin{cor}\label{t:forcingntourn}
If a digraph $D = (V, E)$ satisfies $|E| > |V| \log_2 |V|$ and $D$ is not transitive, it does not have the tournament forcing property.
\end{cor}

\subsection*{Organization}
We begin in~\cref{s:antisid} by studying the tournament anti-Sidorenko property; a primary focus in this section is proving~\cref{c:notanti} and giving several asymptotically tight constructions. Along the way, we observe some connections between our methods and related extremal combinatorics problems and in~\cref{ss:related}, and study a few related questions, including giving bounds on the biclique covering number of a graph.

We proceed in~\cref{s:tsid} to study the tournament Sidorenko property, giving several families of graphs and recursive constructions that yield tournament Sidorenko digraphs (including~\cref{c:growspecial,c:transminusedge}). In~\cref{s:trees}, we study the above pair of Sidorenko properties on directed trees, establishing~\cref{t:star,p:partialtree}, and in~\cref{s:forcing}, we relate these Sidorenko properties to the tournament forcing property. We conclude in~\cref{s:conc} with some open questions and conjectures.

\section{Preliminaries}
We begin by making some observations about the tournament Sidorenko property and the tournament anti-Sidorenko property. For digraphs $D$ and $G$, a labeled copy of $D$ in $G$ is an embedding $\phi: V(D)\to V(G)$ such that if $(x, y)\in E(D)$ then $(\phi(x), \phi(y))\in E(G)$. We will use $N_{L}(D,G)$ to denote the number of labeled copies of $D$ in $G$. 

\begin{prop}\label{p:tanti}
A digraph $D$ is tournament anti-Sidorenko if and only if for every tournament $T$ on $n$ vertices,
$$N_L(D, T) \le 2^{-e(D)} n^{v(D)}.$$
Also, $D$ is tournament Sidorenko if and only if for every tournament $T$ on $n$ vertices,
$$N_L(D, T) \ge (1 - o(1))  2^{-e(D)} n^{v(D)}.$$
\end{prop}
\begin{proof}
We first consider the tournament anti-Sidorenko setting. If $D$ is tournament anti-Sidorenko, i.e. $t_D(T) \le 2^{-e(D)}$, then $N_L(D, T) \le h_D(T) =t_D(T)n^{v(D)}\le 2^{-e(D)} n^{v(D)}$, giving the trivial direction of the implication.

In the other direction, suppose $D$ is not tournament anti-Sidorenko, i.e. there is a tournament $T^*$ such that $t_D(T^*) > 2^{-e(D)}$. Let $k$ denote the number of vertices of $T^*$ and $\epsilon=t_D(T^*) - 2^{-e(D)}$. Let $n$ be a multiple of $k$ such that $n > 2v(D)^2/\epsilon$. Let $T_n$ be an $n$-vertex tournament given by taking a balanced $(n/k)$-blowup of $T^*$ (that is, the lexicographic product of $T^*$ with an empty graph on $n/k$ vertices), and arbitrarily adding a directed edge between each pair of non-adjacent vertices. By construction, $t_D(T_n) \ge t_D(T^*)$ and thus, 
$$N_L(D, T_n) \ge h_D(T_n) - v(D)^2 n^{v(D) - 1} \ge \left(2^{-e(D)} + \eps\right) n^{v(D)} - v(D)^2 n^{v(D) - 1} > \left(2^{-e(D)} + \frac{\eps}{2} \right)n^{v(D)},$$
where the first inequality uses the fraction of maps from $V(D)$ to $V(T_n)$ which are not one-to-one is less than $v(D)^2 /n$. This fact follows from considering the probability that a random map from $V(D)$ to $V(T_n)$ is not one-to-one is at most $\binom{v(D)}{2}$ times the probability that a particular pair of vertices of $D$ get mapped to the same vertex in $V(T_n)$. The above inequality shows the other desired direction. 

Hence, for a digraph $D$, we have $t_D(T) \le 2^{-e(D)}$ for all tournaments $T$ if and only if $N_L(D, T) \le 2^{-e(D)} n^{v(D)}$ for all $n$-vertex tournaments $T$.

Now we prove the claim for the Sidorenko case. Note that $0 \le h_D(T) - N_L(D, T) \le O_D(n^{v(D)-1}) = o(n^{v(D)})$. Since these counts are off by a lower order term (as $n \to \infty$), it follows that $D$ is tournament Sidorenko if and only if $N_L(D,T) \geq (1-o(1))2^{-e(D)}n^{v(D)}$ for every tournament $T$ on $n$ vertices. 
\end{proof}

Thus, we use the alternative equivalent forms of the tournament Sidorenko and anti-Sidorenko properties  from \cref{p:tanti} interchangeably with those stated in~\cref{d:tsidorenko} throughout the remainder of the article.

We also observe that $D$ is tournament Sidorenko (resp. anti-Sidorenko), if and only if the digraph $D'$ obtained by reversing all of the edges of $D$ is. This is because for any tournament $T$, $t_D(T) = t_{D'}(T')$, where $T'$ is the tournament obtained by reversing all of the edges of $T$. We make use of this simple fact throughout the paper.

Finally, we define a strengthening of the tournament anti-Sidorenko property that is useful in proving that certain digraphs constructed from smaller digraphs are tournament anti-Sidorenko.

\begin{defn}[Strongly tournament anti-Sidorenko] \label{d:strong-anti-sid}
Given a digraph $D$ and an independent set $I \subset V(D)$, we say that the pair $(I, D)$ is \textit{strongly tournament anti-Sidorenko} if for any $n$-vertex tournament $T$ and any embedding $\phi: I \hookrightarrow V(T)$, the number of labeled copies $\varphi$ of $D$ in $T$ where $\varphi|_I = \phi$ is at most $2^{-e(D)} n^{v(D) - |I|}.$
\end{defn}

Observe that if $(I, D)$ is strongly tournament anti-Sidorenko, then for any $I' \subseteq I$, $(I', D)$ is also strongly tournament anti-Sidorenko, by an averaging argument over all possible choices of how the vertices in $I\setminus I'$ are embedded. In the special case $I = \emptyset$, $(\emptyset, D)$ being strongly tournament anti-Sidorenko is equivalent to $D$ being tournament anti-Sidorenko in the sense of~\cref{d:tsidorenko}.

Note that the analogous strengthening of the tournament Sidorenko property does not hold for interesting digraphs. In particular, if $(I, D)$ is strongly tournament Sidorenko, then every vertex in $I$ is an isolated vertex in $D$. Indeed, consider a pair $(I, D)$ for which there is a vertex $v\in I$ that is not isolated in $D$. Letting $T$ be a transitive tournament, embedding $v$ as the source or sink, we can an embedding of $I$ which does not extend to any copy of $D$, which shows that the tournament Sidorenko property does not hold for the pair.

\subsection*{Notation}
We let $[n] = \{1, 2, \ldots, n\}$. We use $A \sqcup B$ to denote the disjoint union of two sets $A, B$, and given a directed graph $G = (V, E)$ and some subset $S \subset V$, let $G[S]$ denote the induced subgraph of $G$ on vertex set $S.$ 

Recall that throughout this article, our digraphs (directed graphs) are all oriented, lacking self-loops or antiparallel edges. Given a digraph $D$, we let $\overline{D}$ denote the \textit{underlying undirected graph} of $D$. Throughout the paper, we use standard asymptotic abbreviations $f(n)=O(g(n))$ to mean that there exist $C$ and $n_{0}$ such that $|f(n)|\leq Cg(n)$ for all $n\geq n_{0}$, and $f(n) = o(g(n))$ to mean that $\lim_{n\to\infty} \frac{f(n)}{g(n)}=0$. We will interpret $o(1)$ generally as a quantity that approaches $0$ as $n \to \infty$.

For a graph $G = (V, E),$ and positive integer $k$, we let $V_k$ denote the set of ordered $k$-tuples of distinct vertices $v \in V.$ For $S \in V_k$, we use $S' \subset S$ to denote an \textit{ordered} sub-tuple of the tuple $S$.

Throughout the article, we will use overset or side-by-side labels to indicate the justification for particular inequalities in long chains of reasoning. We also occasionally omit floors and ceilings where they are unimportant for clarity of exposition.

\section{Tournament anti-Sidorenko digraphs}\label{s:antisid}

\subsection{The maximum size of a tournament anti-Sidorenko digraph of a given order}
We first prove~\cref{c:notanti}. It shows that only relatively sparse graphs can be tournament anti-Sidorenko and gives an asymptotically tight edge threshold for this property to possibly hold. The necessity of sparsity follows from the following construction of a family of tournaments where a fixed digraph is systematically over-represented.

\begin{lemma}\label{l:denseoverrep}
Let $D$ be a digraph with $k$ vertices. For each positive integer $n$ that is a multiple of $k$, there is a tournament $T$ on $n$ vertices with $t_D(T) \geq k^{-k}$. 
In particular, if $e(D) \geq k \log_2 k$, then $t_D(T) \geq k^{-k} > 2^{-e(D)}$, and hence $D$ is not tournament anti-Sidorenko. 
\end{lemma}
\begin{proof} 
Let $T$ be a tournament on $n$ vertices given by taking a balanced $(n/k)$-blowup of $D$ (that is, the lexicographic product of $D$ with an empty graph on $n/k$ vertices), and arbitrarily adding a directed edge between non-adjacent vertices. By construction, $t_D(T) \ge \frac{(n/k)^k}{n^k} =k^{-k}$. 

If $e(D) \geq k \log_2 k$, then $t_D(T) \ge k^{-k}=  2^{e(D) - k \log_2 k} \cdot 2^{-e(D)} > 2^{-e(D)}$, as desired. 
\end{proof}

\begin{rem} 
It is not difficult to see that the bound in \cref{l:denseoverrep} is not sharp. For example, we can strengthen the construction given above as follows. Fix digraph $D$ on vertex set $[k]$. By Tur\'an's theorem, $D$ contains an independent set $S \subset [k]$ of size $\alpha(D) \ge k^2/(2e(D) + k)$. 
As in~\cref{l:denseoverrep}, we construct an $n$-tournament as follows, assuming that $k^2 \mid n$. First, take a balanced $(n/k)$-blowup of $D$ to obtain a directed graph on $n$ vertices (assuming $k | n$), with $V(D) = V_1 \sqcup \cdots \sqcup V_k$ with parts corresponding to the vertices of $D$. Note that the induced subgraph on $W = \bigcup_{i \in S} V_i$ is an independent set of size $n \alpha(D) / k$. We can thus equipartition $W = W_1 \sqcup \cdots \sqcup W_k$ into $k$-parts of size $n \alpha(D)/k^2$ and embed a balanced $(n\alpha(D)/k^2)$-blowup of $D$ into these parts. We finally add arbitrarily directed edges between any remaining pairs of non-adjacent vertices to arrive at an $n$-vertex tournament $T$. This modified construction yields a homomorphism density of at least
$$t_D(T) \ge \frac{(n/k)^k}{n^k} + \frac{(n \alpha(D)/k^2)^{k}}{n^{k}} \ge k^{-k} + \left( \frac{\alpha(D)}{k^2} \right)^{k} \ge k^{-k} + \left(2e(D) + k\right)^{-k}.$$
This construction can itself be further improved, but such improvements do not appear to yield an absolute constant factor improvement for any especially large and interesting class of digraphs.
\end{rem}

The above lemma tells us that any digraph $D$ with $e(D) > v(D) \log_2 v(D)$ cannot be tournament anti-Sidorenko. We show that this bound is asymptotically sharp via two very different looking constructions below. We begin with the following simple construction that comes close to the upper bound and also illustrates that tournament anti-Sidorenko digraphs can have any digraph as an induced subgraph.

\begin{lemma}\label{p:anti1}
Let $A$ be an undirected graph with $k$ vertices and $e$ edges and let $D$ be the digraph consisting of the disjoint union of all $2^e$ possible orientations of $A$, so $D$ has $k2^e$ vertices and $e2^e$ edges. Then $D$ is tournament anti-Sidorenko.
\end{lemma}
\begin{proof}
Observe that 
$$\sum_{H : \overline{H} = A} N_L(H, T) = N_L(A, K_n) = n(n-1) \cdots (n-k+1) \le n^k.$$
Furthermore, we see that 
\begin{align*}
N_L(D, T) &\le \prod_{H : \overline{H} = A} N_L(H, T) \overset{\text{AM-GM}}\le \left(2^{-e} \sum_{H \mid  \overline{H} = A} N_L(H, T) \right)^{2^e} \le \left(2^{-e} n^k \right)^{2^e} = 2^{-e(D)}n^{v(D)},
\end{align*}
verifying that $D$ is tournament anti-Sidorenko.
\end{proof}

We next proceed to the two constructions that asymptotically meet our upper bound.
\begin{lemma}\label{p:anti2}
Let $D = (A \sqcup B, F)$ be an orientation of the complete bipartite graph with parts of size $|A| = k, |B| = 2^k$, such that for every subset $A' \subset A$, there is a unique vertex $v \in B$, with $N^+(v) = A'$. Then $(A, D)$ is strongly tournament anti-Sidorenko.
\end{lemma}
\begin{proof}
Note that $v(D)=k+2^k$ and $e(D)=k2^k$. Let $T = (V, E)$ be a tournament on $n$ vertices. Let $S$ be a fixed ordered $k$-tuple of vertices of $V$. Then, the number of labelled copies of $D$ in $T$ where the image of $A$ is $S$ is given by
\begin{align*}
& \prod_{S' \subset S} \left|\{v \in V \backslash S \mid N^+(v) \cap S = S'\}\right| \\
&\overset{\text{AM-GM}}\le \left(2^{-k} \sum_{S' \subset S} \left|\{v \in V \backslash S \mid N^+(v) \cap S = S'\}\right| \right)^{2^k} \\
&= \left(2^{-k} | V \backslash S| \right)^{2^k} \le \left(2^{-k} n\right)^{2^k} = 2^{-e(D)}n^{v(D) - k}.
\end{align*}
Since this holds for every choice of $S$, we obtain the desired result.
\end{proof}

Our final construction of ``relatively dense'' tournament anti-Sidorenko digraphs is a special case of the following more general method to construct larger tournament anti-Sidorenko digraphs from smaller instances. 

\begin{defn}
Given a digraph $G$, a labeled digraph $D$ on $k$ vertices and an ordered $k$-tuple $S$ of vertices of $G$, we use $G[S] \sim D$ to mean that the induced sub-digraph $G[S]$ contains a \textit{labeled} copy of $D$ that respects the given ordering of vertices in $S$.
\end{defn}

\begin{prop}\label{t:antiextend}
Let $D_1, D_2$ be a pair of tournament anti-Sidorenko digraphs such that there is an independent subset $I \subset V(D_1)$ with $|I| = v(D_2)$. If all vertices in $I$ each have the same in- and out-neighborhood, then the digraph $D$ obtained by embedding $D_2 \hookrightarrow I$ is also tournament-anti-Sidorenko.
\end{prop}
\begin{proof}
Let $|I| = k = v(D_2)$ and let $D_3 = D_1[V(D_1) \backslash I]$ be the induced subgraph of $D_1$ on the vertex subset excluding $I$. We further define $D_4 = D_1[V(D_3) \cup \{v\}]$ to be the induced subgraph of $D_1$ on vertex set $V(D_3)$ along with an arbitrarily chosen additional vertex $v \in I$ (any choice of $v \in I$ yields the same digraph by the condition on $I$ in the statement).

Let $T$ be an $n$-vertex tournament.
We count labeled copies of $D_1$ in $T$ as follows: we identify an ordered $v(D_3)$-tuple, $S$, of vertices of $T$, such that $T[S] \sim D_3$ and let $U_S$ be the set of vertices $u \in V(T)$ such that $T[(S, u)] \sim D_4$, i.e. $u \in U_S$ if when taken together with $S$, contains a labeled copy of $D_4$ where the image of $\{v\} \in D_4$ is $u$. With these definitions, we observe that
$$h_{D_1}(T) \ge \sum_{S \mid T[S] \sim D_3} (|U_S|)^k,$$
since any sequence of $k$ vertices from $U_S$ together with $S$ gives a homomorphic copy of $D_1.$
We are able to count copies of $D$ in $T$ and show that $D$ is tournament anti-Sidorenko as follows.
\begin{align*}
N_L(D, T) &= \sum_{S \mid T[S] \sim D_3} N_L(D_2, T[U_S]) \\
&\le 2^{-e(D_2)} \sum_{S \mid T[S] \sim D_3} (|U_S|)^{k} \tag{$D_2$ t.-anti.-Sid.} \\
&\le 2^{-e(D_2)} h_{D_1}(T) \\
&\le 2^{-e(D)} n^{v(D)}. \tag{$D_1$ t.-anti.-Sid.}
\end{align*}
\end{proof}
With this result in hand, we can complete our third construction.

\begin{defn}
Setting $S_{0} = \emptyset$, we recursively define the \textit{iterated balanced star} $S = S_{k + 1}$ as a graph with $k+1$ vertices $V(S) = \{v\} \cup U \cup W$ with $|U| = \lfloor k/2\rfloor, |W| = \lceil k/2 \rceil$, with edge set consisting of a balanced star with central vertex $v$, in-edges to $v$ coming from $U$, out-edges from $v$ going to $W$, and let $S_{k+1}[U] = S_{|U|}, S_{k+1}[W] = S_{|W|}$.
\end{defn}
\begin{rem}\label{rem:balanced}
    Note that for $k\ge 3$, $e(S_k)$ satisfies the following recursion:
    \[e(S_k) = k-1 + e(S_{\lfloor \frac{k-1}{2}\rfloor}) + e(S_{\lceil \frac{k-1}{2}\rceil}).\]
    We now prove by induction that $e(S_k) = t(k+1) - 2^{t+1} + 2$ where $t$ is the largest integer with $2^t \le k+1$. This is true for $e(S_1) = 0$ and $e(S_2) = 1$. For $k\ge 3$, note that if $t$ is the largest integer with $2^t \le k+1$, then we see that $2^{t-1} \le \lfloor \frac{k-1}{2}\rfloor + 1 \le \lceil \frac{k-1}{2}\rceil + 1 < 2^t$. Hence, by induction hypothesis we have $e(S_{\lfloor \frac{k-1}{2}\rfloor}) = (t-1)(\lfloor \frac{k-1}{2}\rfloor + 1) - 2^t + 2$ and $e(S_{\lceil \frac{k-1}{2}\rceil}) = (t-1)(\lceil \frac{k-1}{2}\rceil + 1)-2^t + 2$. Thus,
    \[e(S_k) = k-1 + (t-1)(k+1) - 2\cdot 2^{t} + 4 = t(k+1) - 2^{t+1} + 2.\]
    In particular, $\log_2(k+1) - 1 < t \le \log_2(k+1)$ implies that $e(S_k) \ge (k+1)\log_2(k+1) - 2k$.
\end{rem}

\begin{lemma}\label{p:anti3}
The iterated balanced star $S_{k+1}$ is tournament anti-Sidorenko.
\end{lemma}
\begin{proof}
We proceed inductively.
Note that $S_3$ is an orientation of $K_{1,2}$ with a central vertex that has both in-degree one and out-degree $1$. Therefore, by~\cref{t:star} (which we verify later), this digraph is tournament anti-Sidorenko. For any $k > 2$, we assume that $S_i$ is anti-Sidorenko for each $i \le k$. Now, let $D$ be an orientation of the $k+1$ vertex star that is tournament anti-Sidorenko per~\cref{t:star} with central vertex $v \in V(D)$ that has (without loss of generality) in-degree $\lfloor k /2 \rfloor$ and out-degree $\lceil k/2 \rceil$. We let $I = N^-(v) \subset V(D)$, and apply~\cref{t:antiextend} with $D_1 = D, D_2 = S_{\lfloor k/2 \rfloor}$ to obtain a resulting tournament anti-Sidorenko digraph $D'$. We then apply~\cref{t:antiextend} again, this time letting $I = N^+(v) \subset V(D)$, letting $D_1 = D', D_2 = S_{\lceil k/2 \rceil}$. The resulting digraph, which is $S_{k+1}$, is consequently tournament anti-Sidorenko.
\end{proof}

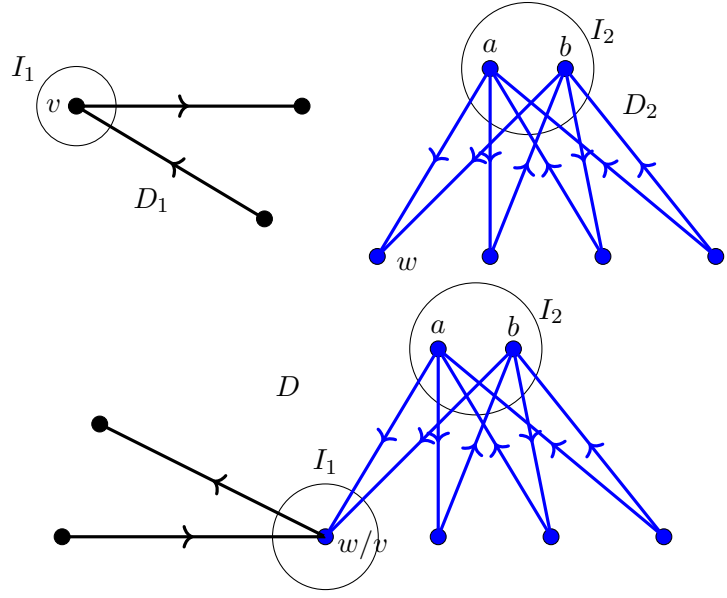
\begin{figure}[t]
\centering
\begin{tikzpicture}
\draw[fill=black] (4.5,2) circle (3pt);
\draw (4.5, 2) ellipse (15pt);
\draw[fill=black] (7.5,2) circle (3pt);
\draw[fill=black] (7,0.5) circle (3pt);
\draw[fill=blue]  (10,2.5) circle (3pt);
\draw[fill=blue]  (11,2.5) circle (3pt);
\draw (10.5, 2.5) ellipse (25pt);
\draw[fill=blue] (8.5, 0) circle (3pt);
\draw[fill=blue] (10, 0) circle (3pt);
\draw[fill=blue] (11.5, 0) circle (3pt);
\draw[fill=blue]  (13, 0) circle (3pt);
\node at (3.8,2.5) {$I_1$};
\node at (11.5,3) {$I_2$};
\node at (10,2.8) {$a$};
\node at (11,2.8) {$b$};
\node at (4.2, 2) {$v$};
\node  at (8.9,-0.1) {$w$};
\node at (5.5, 0.75) {$D_1$};
\node at (12, 2) {$D_2$};
\begin{scope}[very thick,decoration={
    markings,
    mark=at position 0.5 with {\arrow{>}}}
    ] 
     \draw[black,postaction={decorate}] (4.5,2)--(7.5, 2);
     \draw[black,postaction={decorate}] (7,0.5)--(4.5, 2);
     \draw[blue,postaction={decorate}] (10,2.5)--(8.5, 0);
     \draw[blue,postaction={decorate}] (11,2.5)--(8.5, 0);
    \draw[blue,postaction={decorate}] (10,2.5)--(10, 0);
     \draw[blue,postaction={decorate}] (10, 0)--(11,2.5);
        \draw[blue,postaction={decorate}] (11,2.5)--(11.5, 0);
     \draw[blue,postaction={decorate}] (11.5,0)--(10,2.5);
    \draw[blue,postaction={decorate}] (13,0)--(11, 2.5);
    \draw[blue,postaction={decorate}] (13,0)--(10, 2.5);
\end{scope}
\end{tikzpicture}
\begin{tikzpicture}
\draw[fill=black] (5.5,1.5) circle (3pt);
\draw[fill=black] (5,0) circle (3pt);
\draw[fill=blue]  (10,2.5) circle (3pt);
\draw[fill=blue]  (11,2.5) circle (3pt);
\draw (8.5,0) ellipse (20pt);
\draw (10.5, 2.5) ellipse (25pt);
\draw[fill=blue] (8.5, 0) circle (3pt);
\draw[fill=blue] (10, 0) circle (3pt);
\draw[fill=blue] (11.5, 0) circle (3pt);
\draw[fill=blue]  (13, 0) circle (3pt);
\node at (8, 2) {$D$};
\node at (8.5,1.0) {$I_1$};
\node at (11.5,3) {$I_2$};
\node at (10,2.8) {$a$};
\node at (11,2.8) {$b$};
\node  at (9,-0.1) {$w / v$};
\begin{scope}[very thick,decoration={
    markings,
    mark=at position 0.5 with {\arrow{>}}}
    ] 
    \draw[black,postaction={decorate}] (5,0)--(8.5, 0);
    \draw[black,postaction={decorate}] (8.5,0)--(5.5, 1.5);
     \draw[blue,postaction={decorate}] (10,2.5)--(8.5, 0);
     \draw[blue,postaction={decorate}] (11,2.5)--(8.5, 0);
    \draw[blue,postaction={decorate}] (10,2.5)--(10, 0);
     \draw[blue,postaction={decorate}] (10, 0)--(11,2.5);
        \draw[blue,postaction={decorate}] (11,2.5)--(11.5, 0);
     \draw[blue,postaction={decorate}] (11.5,0)--(10,2.5);
    \draw[blue,postaction={decorate}] (13,0)--(11, 2.5);
    \draw[blue,postaction={decorate}] (13,0)--(10, 2.5);
\end{scope}
\end{tikzpicture}
\caption{Per~\cref{p:even-star}, we have that $(I_1, D_1)$ is strongly tournament anti-Sidorenko with $I_1 = \{w\}$, and per~\cref{p:anti2} $(I_2, D_2)$ is strongly tournament anti-Sidorenko with $I_2 = \{a, b\}$. Then by~\cref{l:glue-anti-sid} the digraph $D$ formed by identifying $v \in D_1$ with $w \in D_2$ satisfies $(I_2, D)$ is strongly tournament anti-Sidorenko.}
\label{f:growstrong}
\end{figure}

\begin{proof}[Proof of~\cref{c:notanti}]
\cref{t:antiextend} implies $f(k) \le k \log_2 k$.
By considering the construction of~\cref{p:anti3} (in conjunction with the observation~\cref{rem:balanced}) or the construction~\cref{p:anti2}, we obtain the lower bound $f(k) \ge k\log_2k - O(k)$. 
\end{proof}
We note some additional digraphs we can show to have the tournament anti-Sidorenko property via similar methods to above.
Below we construct larger digraphs that have the (strong) tournament anti-Sidorenko property from smaller digraphs. Specifically, given two digraphs $D_1, D_2$ and an independent set $I_1 \subseteq V(D_1)$, we say that digraph $D$ is formed by \textit{identifying} $|I_1|$ vertices of $D_2$ with $I_1 \subseteq V(D_1)$, if, with slight abuse of notation, we can view vertices in $D_1$ and $D_2$ as vertices in $D$ such that $V(D) = V(D_1) \cup V(D_2)$, $V(D_1) \cap V(D_2) = I_1$ and $E(D) = E(D_1) \sqcup E(D_2)$. We give an example in~\cref{f:growstrong} of this identification process and the resulting strongly tournament anti-Sidorenko digraph exhibited by~\cref{l:glue-anti-sid}.
\color{black}

\begin{lemma}\label{l:glue-anti-sid}
If both $(I_1, D_1)$ and $(I_2, D_2)$ are strongly tournament anti-Sidorenko pairs (per~\cref{d:strong-anti-sid}) with $v(D_2) \ge |I_1|$, then the digraph $D$ formed by identifying $|I_1|$ vertices of $D_2$ with $I_1 \subset V(D_1)$ satisfies that $(I_2, D)$ is strongly tournament anti-Sidorenko.
\end{lemma}
\begin{proof}
As mentioned above, we denote vertices in $D$ by vertices in $D_1$ and $D_2$. In particular, we may assume that $V(D) = V(D_1)\cup V(D_2)$, $I_1$ and $I_2$ are both subsets of $V(D)$, and $V(D_1)\cap V(D_2) = I_1$.

Let $T$ be an $n$-vertex tournament with vertex set $V$. Given any embedding $\phi: I_2 \hookrightarrow V$, we count the number of labeled copies $\varphi: V(D) \to V$ of $D$ in $T$ such that $\varphi|_{I_2} = \phi$.

As $(I_2, D_2)$ is strongly tournament anti-Sidorenko, the number of labeled copies $\varphi' = \varphi|_{V(D_2)}: V(D_2)\hookrightarrow V$ of $D_2$ in $T$ such that $\varphi|_{I_2} = \phi$ is at most $2^{-e(D_2)}n^{v(D_2) - |I_2|}$.

For each $\varphi'$, as $(I_1, D_1)$ is strongly tournament anti-Sidorenko, the number of labeled copies $\varphi'' = \varphi|_{V(D_1)}: V(D_1)\hookrightarrow V$ of $D_1$ in $T$ (with vertices in $\varphi'(V(D_2)\setminus I_1)$ removed) such that $\varphi''|_{I_1} = \varphi'|_{I_1}$ is at most $2^{-e(D_1)}(n - |V(D_2)\setminus I_1|)^{v(D_1)-|I_1|} \le 2^{-e(D_1)}n^{v(D_1)-|I_1|}$.

Note that since $V(D) = V(D_1)\cup V(D_2)$, $(\varphi', \varphi'')$ uniquely determines $\varphi$. Hence, the total number of ways to choose $\varphi$, which is the number of ways to choose the pair $(\varphi', \varphi'')$, is at most 
\[2^{-e(D_1) - e(D_2)} n^{v(D_1) - |I_1| + v(D_2) - |I_2|} = 2^{-e(D)} n^{v(D) - |I_2|}.\]
This establishes that $(I_2, D)$ is strongly tournament anti-Sidorenko.
\end{proof}

\begin{cor} \label{p:tight-anti-sid2}
    Consider digraph $D = (A_1 \sqcup A_2, F)$ given in~\cref{p:anti2} with parts of size $|A_1| = s, |A_2| = 2^s$. Let $B$ be a tournament anti-Sidorenko digraph with $s$ vertices. Then for any choice of embedding $V(B) \hookrightarrow V(A_1)$, the digraph $D' = (A_1 \sqcup A_2, F \sqcup E(B))$ is also tournament anti-Sidorenko.
\end{cor}
\begin{proof}
This follows by applying~\cref{l:glue-anti-sid} with $I_1 = A_1, D_1 = D$, $I_2 = \emptyset$ and $D_2 = B$, noting that $|B| = |A_1|$ and thus we identify each vertex in $B$ with a distinct vertex in $A_1 \subset V(D_1).$
\end{proof}

Note that by choosing $B$ to be an tournament anti-Sidorenko digraph on $s$ vertices with a maximum number of edges, the construction in~\cref{p:tight-anti-sid2} has $k = s + 2^s$ vertices and $f(s) + s 2^s = k \log k - O(k)$ edges for $f(\cdot)$ as defined in~\cref{c:notanti}, and thus is a tight example for~\cref{c:notanti} up to the $O(k)$ term. The construction in~\cref{p:anti3} is also a tight example for~\cref{c:notanti}.

\begin{prop}
    Let $D$ be a tournament anti-Sidorenko digraph. If there exists a pair of vertices $v,w \in V(D)$ such that $(v,w),(w,v)\notin E(D)$ and such that the digraphs $D'=(V(D),E(D)\cup (w,v))$ and $D'' = (V(D),E(D)\cup (v,w))$ are isomorphic, $D'\cong D''$, then $D'$ (equivalently, $D''$) is also tournament anti-Sidorenko.
\end{prop}
\begin{proof}
Let $T$ be an $n$-vertex tournament.
Note that if $D' \cong D''$, then $N_L(D', T) = N_L(D'', T)$. Further note that $N_L(D, T) = N_L(D', T) + N_L(D'', T)$. Therefore, if $N_L(D, T) \le 2^{-e(D)} n^{v(D)},$ then $N_L(D', T) = N_L(D'', T) = \frac12 N_L(D, T) \le 2^{-e(D) - 1} n^{v(D)},$ showing that $D', D''$ are tournament anti-Sidorenko.
\end{proof} This proposition also holds, by the same argument, when we replace anti-Sidorenko by Sidorenko. We obtain the following direct corollary (noting that $C_{2k}$ is tournament anti-Sidorenko for $k$ odd per~\cite{GKLV20}). 
\begin{cor}
Consider a directed cycle with $2k$ vertices $C_{2k}$ where $k$ is odd, with vertices $v_1, \ldots, v_{2k}$, Let $D = (\{v_1, \ldots, v_{2k}\}, E(C_{2k}) \cup \{(v_1, v_{k+1})\}),$ (i.e. we add a directed edge that connects a pair of vertices at maximum distance). The resulting digraph is tournament anti-Sidorenko.
\end{cor}

\subsection{Related questions}\label{ss:related}
In the proof that any tournament anti-Sidorenko digraph $D$ has at most $v(D) \log_2 v(D)$ edges (\cref{l:denseoverrep}), we constructed a family of tournaments $F_k$ from balanced blowups of $D$ such that each tournament $F_k$ has significantly more copies of $D$ than expected in a randomly oriented tournament on the same vertices. In fact, we could have instead chosen to construct a family of tournaments $(F_n)_{n = 1}^{\infty}$ from any digraph $D'$ such that there exists a homomorphism $D \to D'$. If $k' := v(D') < v(D) =: k$, this yields an improved upper bound, as the argument of~\cref{l:denseoverrep} implies that $$t_D(F_n') \ge n^k/(k' n)^k = 2^{-k \log_2 k'} > 2^{-e(D)},$$ whenever $e(D) > k \log_2 k'$, that $D$ cannot be tournament anti-Sidorenko if it has at least $v(D) \log_2 v(D')$ edges. We also obtain an improved bound if $D$ has \textit{more than one homomorphism} into some tournament on $k$ vertices, as we can choose to blow up such a tournament instead of an arbitrarily ``filled-in'' copy of $D$.

Thus, the study of the ``densest'' possible tournament anti-Sidorenko digraphs is intimately related to the following more general question:
\begin{qn}
What is the minimum number of edges in a digraph with $k$ vertices that has at most one homomorphism into any tournament on $k$ vertices?
\end{qn}
This question is natural in its own right, but also has close connections to other questions in extremal combinatorics, like the size of a biclique covering. We investigate some of these connections in the remainder of this subsection.
We begin by studying a slightly weaker question, where we replace homomorphism above with embedding.
\begin{lemma}
    Let $D$ be a digraph on $k$ vertices. Suppose that $D$ has at most one embedding into any $k$-vertex tournament. Then $e(D) \ge k\log_2 k - O(k)$.
\end{lemma}
\begin{proof}
    Consider a uniform random tournament $T$ on $[k]$, and let $\pi$ be any permutation $V(D)\to [k]$. The probability that $\pi$ is an embedding is exactly $1/2^{e(D)}$. Since $D$ has at most one embedding, we have that $k!/2^{e(D)}\le 1$. This gives $e(D)\ge \log_2 k! = k\log_2 k - O(k)$.
\end{proof}

We proceed to give a lower bound on the total number of vertices of a biclique cover of a graph. 
We let 
$K_{A, B}$ be the \textit{undirected} complete bipartite graph with vertex set $A \sqcup B$ and edge set $\{(a, b) \mid a \in A, b \in B\}$.
\begin{defn}
Given an undirected graph $H = ([n], E)$, a \textit{biclique cover} of $H$ is a collection of complete bipartite graphs $\{K_{A_i, B_i}\mid A_i, B_i \subset [n],\, A_i \cap B_i = \emptyset\}_{i \in [k]}$ that are each subgraphs of $H$ such that for each $e \in E$, there is an $i \in [k]$ such that $e \in E(K_{A_i, B_i})$.
\end{defn}

\begin{lemma}\label{lem:biclique-cover}
Let $H = ([n], E)$ be an undirected graph with $|E| \ge \binom{n}{2} - s$ where $s$ is a nonnegative integer. Let $\{K_{A_i, B_i}\mid A_i, B_i \subset [n],\, A_i \cap B_i = \emptyset\}_{i \in [k]}$ be a biclique cover of $H$. Then 
$$\sum_{i=1}^k |A_i| + |B_i|\ge n\log_2n - n\log_2\frac{s+n}{n} - O(n).$$
\end{lemma}
\begin{proof}
Let $C_i = A_i \sqcup B_i$ for $i \in [k]$. For each vertex $v \in [n] = V(H),$ let $f(v)$ be the number of bicliques that include vertex $v$, so $f(v) = |\{i \in [k] \mid v \in C_i\}|$. Now let $T_t = \{v \in [n] \mid f(v) = t\}$ be the collection of vertices in $H$ that are covered by exactly $t$ bicliques in the cover.
We first observe that each $|T_t|$ is not too big.
\begin{claim}\label{c:tt-small}
For each $t \in [k]$, $|T_t| \le 2^t + \sqrt{s 2^{t+1}}$.
\end{claim}
\begin{proof}
Let $X_i$ be either $A_i, B_i$ chosen uniformly at random independently for $i\in [k]$. Let $Y = T_t\setminus (X_1\cup \cdots \cup X_k)$. For each $v\in T_t$, since $v$ is in exactly $t$ sets $C_i$, $v$ is included in $Y$ with probability $2^{-t}$. Therefore $\mathbb E[|Y|] = |T_t|/2^t$. Let $H^c$ be the complement of $H$, so that $e(H^c) \le s$. Now we claim that $H^c[Y]$ is the complete graph on $Y$. For any two distinct vertices $u, v\in Y$ we know that there is no such $i$ such that $u\in A_i$ and $v\in B_i$ (otherwise exactly one of them is excluded by $X_i$), so $\{u, v\}\notin E(K_{A_i, B_i})$ for all $i\in [k]$. Since the $K_{A_i, B_i}$ form an edge clique cover, $\{u, v\}\in E(H^c)$; this proves $H^c[Y]$ is the complete graph on $Y$. Thus, we have that
    \[\mathbb{E}[e(H^c[Y])] = \mathbb{E}\left[\binom{|Y|}{2}\right] = \frac12\mathbb{E}[|Y|^2] - \frac12\mathbb{E}[|Y|]
    \overset{\text{Cauchy-Schwarz}}\ge 
    \frac12\mathbb{E}[|Y|]^2 - \frac12\mathbb{E}[|Y|]=
    \frac{|T_t|^2}{2^{2t+1}} - \frac{|T_t|}{2^{t+1}}.\]
    On the other hand, for each edge $\{u, v\}\in H^c[T_t]$, the probability that it is an edge of $H^c[Y]$ is $\Pr[u\in Y\text{ and } v\in Y] \le \Pr[u\in Y] = 2^{-t}$. Therefore, we conclude that 
    $$
    \frac{s}{2^t} \ge \frac{e(H^c[T_t])}{2^t} \ge \mathbb{E}[e(H^c[Y])]\ge \frac{|T_t|^2-|T_t|}{2^{2t+1}}.$$
    Rearranging gives $\frac{|T_t|}{2^t} \le \frac{1}{2}\left(1+\sqrt{1+8\frac{s}{2^t}}\right)$. Thus, $|T_t| \le 2^{t-1} + \sqrt{2^{2t-2} + s 2^{t+1}} \le 2^t + \sqrt{s 2^{t+1}}$, as desired.
\end{proof}

We can then bound the number of vertices in our biclique cover. Let $t_0 = \lfloor{2\log_2 n  - \log_2 (s+n)}\rfloor$. Then, we have the following: 
\begin{align*}
     \sum_{i=1}^k |A_i|+|B_i|
     &= \sum_{i=1}^k |C_i| = \sum_{t=1}^k t|T_t| \\
    &\ge t_0n + \sum_{t=1}^{t_0-1}(t-t_0)|T_t| \tag{$\sum_{t=1}^k |T_t| = n$}\\
    & \ge t_0n + \sum_{t=1}^{t_0-1}(t-t_0)\left(\sqrt{s2^{t+1}} + 2^t\right) \tag{\cref{c:tt-small}}\\
    & = t_0n - \sqrt{s2^{t_0+1}}\sum_{\ell=1}^{t_0-1}\ell(\sqrt{2})^{-\ell} - 2^{t_0}\sum_{\ell=1}^{t_0-1}\ell2^{-\ell} \tag{$\ell = t_0-t$}\\
    &\ge t_0n - O(n),
\end{align*}
where in the last inequality we used that by the choice of $t_0$, $2^{t_0} \le \frac{n^2}{s+n} \le n$ and $\sqrt{s2^{t_0}}\le n$ and that the arithmetico-geometric series partial sums are bounded above by absolute constants. The result then follows since $t_0 \ge \log_2 n - \log_2 \left(\frac{s}{n} + 1\right) - 1$.
\end{proof}
\begin{rem}
We remark that the bound in~\cref{lem:biclique-cover} is tight apart from the $O(n)$ term. Consider a graph with $n = 2^r$ vertices, each vertex represented by a binary string of length $r$. Let $k \le r$ be a positive integer. For each $i\in [k]$, let $A_i$ ($B_i$, respectively) be the set of vertices where the $i$-th digit is $0$ ($1$, respectively). Note that $E(K_{A_i, B_i})$ comprises the set of pairs of vertices that differ on the $i$-th bit. Therefore, $\sum_{i=1}^k |A_i|+|B_i| = kn$.

Meanwhile, we see that for the graph $H$ with edge set $E(H) = \bigcup_{i\in [k]}E(K_{A_i, B_i})$, the edges include all pairs of vertices that differ on at least one of the first $k$ bits. Thus, its complement $H^c$ is a disjoint union of $2^k$ complete graphs on $2^{r-k}$ vertices. Letting $s = 2^k\binom{2^{r-k}}{2} = n\frac{2^{r-k}-1}{2}$, so that $r-k = \log_2\frac{n+2s}{n}$, we have $\sum_{i=1}^k |A_i|+|B_i| = n(r-(r-k)) = n\log_2 n - n\log_2\frac{n+2s}{n}$. This matches the bound in~\cref{lem:biclique-cover}.
\end{rem}
\begin{prop}
Let $D$ be a digraph on $k$ vertices. If $D$ has at most one homomorphism to any tournament on $k$ vertices, then $e(D) = k\log_2 k - O(k\log\log k)$.
\end{prop}
\begin{proof}
We may assume that $e(D) < k \log_2 k.$
Let $d = (\log_2 k)^2$. Let $W = \{v_1, \ldots, v_{\ell}\}\subseteq V(D)$ be the collection of vertices with degree more than $d$. Let $X = V(D) \backslash W$, so that all vertices in $D[X]$ have degree at most $d$ and therefore, there are at most $(k - \ell)d^2$ directed paths of length $2$ in $D[X].$
We make use of the following claim:

\begin{claim}\label{lem:one-homo-cond}
Let $D$ be a digraph on $b$ vertices. Suppose that $D$ has at most one homomorphism to any tournament on $k$ vertices. Then, for any distinct $u, v \in V(D)$ either there is an edge on $\{u, v\} \in E(D)$ or there exists a $2$-edge directed path connecting $u$ and $v$.
\end{claim}
\begin{proof}
Suppose to the contrary that there was a pair of vertices $u, v\in V(D)$ that satisfied neither of the above condition, but $D$ had only one homomorphism to any $k$ vertex tournament. We construct a tournament on $k$ vertices as follows. We begin by embedding a copy of $D$. We then add edges $(u, w)$ for all $w \in N^+(v)$ and all edges $(w, u)$ for $w \in N^-(v)$; for all remaining vertex pairs that do not support an edge, we add an arbitrary edge. This yields a $k$-vertex tournament $T$, with at least one homomorphism $D \to T$ given by the pre-fixed embedding $\phi: V(D) \hookrightarrow V(T)$. However, if we consider the map $\phi': V(D) \rightarrow T$ where $\phi'(v) = \phi(u)$ and $\phi'(w) = \phi(w)$ for all other $w \in V(T)$, our choice of $u, v$ means that $N^+(u) \cap N^-(v) = N^-(u) \cap N^+(v) = \emptyset$ and that there is no edge on $\{u, v\}$; thus $\phi'$ is a homomorphism, contradicting our assumption on $D$.
\end{proof}

For each $i \in [\ell],$ let $A_i = N^+(v_i)\cap X$ and $B_i = N^-(v_i)\cap X$. By \cref{lem:one-homo-cond}, we know that any two vertices in $X$ are either adjacent or are connected by a length $2$ directed path. Among all the pairs of vertices in $X$, at most $e(D[X])\le e(D)$ pairs are adjacent; at most $(k-\ell)d^2$ pairs are connected by a directed path of length two, where the intermediate vertex also lies in $X$. Therefore, the remaining pairs of vertices in $X$ must be connected by directed paths through $v_i$ of length $2$ (i.e. with one endpoint in $A_i$ and the other in $B_i$). We can then apply~\cref{lem:biclique-cover} to the auxiliary graph with edge set $\{E(K_{A_i, B_i})\}_{i \in [\ell]}$, and see that for $s = e(D) + (k-\ell)d^2$ we have 
$$\sum_{i = 1}^{\ell} |A_i| + |B_i| \ge |X| \log_2 |X| - |X| \log_2 \frac{s + |X|}{|X|} - O(|X|).$$
Note that $|X| = k - \ell$ and since $d = (\log_2 k)^2,$ we have $k \le 2e(D)/d \le \frac{2k}{\log_2 k} \le \frac{k}{2}$ for $k$ sufficiently large. Further, $s \le e(D) + kd^2 \le 2k(\log_2 k)^4$. By our choice of $A_i, B_i$, $\sum_{i = 1}^{\ell} |A_i| + |B_i| \le e(D[X, W]) \le e(D)$. Combining these inequalities gives
\begin{align*}
     e(B) &\ge \sum_{i=1}^{\ell} |A_i|+|B_i|\\
        & \ge (k-\ell) \log_2 (k/2) - n\log_2 \left(\frac{2s}{k}+1\right) - O(k)\\
        & \ge \left(k-\frac{2k}{\log_2 k}\right)(\log_2k-1) - k\log_2(4(\log_2k)^4+1) - O(k) \\
        & \ge k\log_2 k - O(k\log\log k).\qedhere
\end{align*}
\end{proof}

Finally, we show that our result is essentially tight, using a construction motivated by~\cref{p:anti2}.
\begin{lemma}
For all sufficiently large positive integers $k$, there exists a digraph $D$ on $k$ vertices with $k\log_2 k + O(k)$ edges such that any tournament on $k$ vertices has at most one homomorphic copy of $D$.
\end{lemma}
\begin{proof}
We assume that $k$ is sufficiently large. Let $s = \lceil \log_2 k\rceil$. Let $D$ be a digraph on $k$ vertices with $V(D) = A_1 \sqcup A_2 \sqcup \{x\}$ with $|A_1| = s, |A_2| = t = k - s - 1$. By our choice of $s$, we have $t \le 2^s$. We let $E(D)$ comprise the following edges:
\begin{enumerate}
    \item a transitive tournament on $A_1$;
    \item $\{ (x, u) : u \in A_2\}$;
    \item $\{(v, x) : v \in A_1\}$;
    \item An orientation of the complete bipartite graph with parts $A_1, A_2$ such that: the sets $N^+(u) \cap A_1$ for $u \in A_2$ are distinct and $||N^+(v) \cap A_2| - |N^+(v) \cap A_1|| \le 1$ for all $v \in A_1$.
\end{enumerate}
We show that such a digraph $D$ exists. Equivalently, we show that there exists an orientation of the complete bipartite graph on $A_1 \sqcup A_2$ that satisfies (4). Consider all $2^s$ subsets of $A_1$, $\{S_i\}_{i \in [2^s]}$, ordered so that $S_{2i} = A_1 \backslash S_{2i - 1}$. Let $A_2 = \{u_i\}_{i \in [t]}$. Now we orient edges in $K_{A_1, A_2}$ as follows. We let $N^+(u_i) \cap A_1 = S_i$ and $N^-(u_i) \cap A_1 = A_1 \backslash S_i$. Thus, the sets $N^+(u) \cap A_1$ with $u \in A_2$ are distinct. Each $v \in A_1$ is in exactly one of $S_{2i-1}, S_{2i}$ for each $i$ and thus $v$ is an in-neighbor of exactly one of $u_{2i-1}, u_{2i},$ and an out-neighbor of the other. Thus, when $t$ is even, $|N^+(v) \cap A_2| = |N^+(v) \cap A_1|$, and when $t$ is odd, the neighborhood sizes differ by at most $1.$

Notice that as constructed, the number of edges in (1), (2), (3), (4) are $\binom{s}{2}, t, s, st$ respectively, so $e(D) = \binom{s}{2} + t + s + st = k\log_2 k + O(k)$. We show that any tournament $T$ on $k$ vertices has at most one homomorphic copy of $D$ provided $t > 2s+1$. First, note that any such homomorphism $\phi$ must be injective as if $\phi(u) = \phi(v)$ for $u, v \in V(D)$, we cannot have an edge between $u, v$. Since all vertices in $A_1 \sqcup \{x\}$ are complete to the balance of the graph, this implies $u, v \in A_2$. However, by construction, we can then find $w \in A_1$ with $w \in N^+(u) \backslash N^+(v)$; thus the edge between $\phi(w)$ and $\phi(u) = \phi(v)$ cannot be compatible with both orientations $(u, w), (w, v).$

Thus, any such $\phi$ is injective and satisfies that $N^-_T(\phi(x)) = \phi(A_1)$ and $N^+_T(\phi(x)) = \phi(A_2)$. Similarly, for all $v \in A_1,$ $$|N^-_T(\phi(v))| = |N^-_D(v)| \ge |N^-_D(v) \cap A_2| \ge \frac{t-1}{2}.$$
Suppose to the contrary that there were at least two embeddings $\phi_1, \phi_2: D \hookrightarrow T$. Let $\phi_1(x) = x_1, \phi_2(x) = x_2$ with $x_1 \neq x_2$ such that without loss of generality $(x_1, x_2) \in E(T)$. Because $\phi_1$ is an embedding and $|D| = k = |T|$, and $x$ is adjacent to all other vertices in $D$, the in-degree of $x$ in $D$ is the same as the in-degree of $x_1$ in $T$. With a similar relation for $x_2$, we know that $d^-(x_1) = s = d^-(x_2)$. Therefore, as $(x_1, x_2) \in E(T)$ we have $x_1 \in \phi_2(A_1)$ and thus per above, $d^-(x_1) \ge (t-1)/2$, which yields a contradiction if $t > 2s + 1$. Thus $\phi_1(x) = \phi_2(x)$.

This implies that $\phi_1(A_1) = N^-(\phi_1(x)) = N^-(\phi_2(x)) = \phi_2(A_1)$. As $D[A_1]$ is a transitive tournament, there is at most one way to embed $A_1 \hookrightarrow \phi_1(A_1)$ as a homomorphism. Hence $\phi_1|_{A_1} = \phi_2|_{A_1}$. Finally, if $\phi_1\ne \phi_2$, since they must agree on $A_1\cup \{x\}$, and they are bijections between $V(D)$ and $V(T)$, there must be $v_1\ne v_2\in A_2$ such that $\phi_1(v_1) = \phi_2(v_2)$. We may assume that there exists $u\in A_1$ such that $u\in N^+(v_1)$ and $u\notin N^+(v_2)$ by the distinctness of the out-neighborhoods of $v_i$ in $A_2$. Then, $\phi_1(u) = \phi_2(u)$ and $\phi_1(v_1) = \phi_2(v_2)$, so the edge between them cannot be compatible with both $(v_1, u)$ and $(u, v_2)$. This gives a contradiction, showing $\phi_1 = \phi_2$. Hence, there is at most one homomorphism from $D$ to $T$ for any tournament $T$. \qedhere
\end{proof}
\color{black}

\section{Tournament Sidorenko digraphs}\label{s:tsid}
We next turn our attention to studying the tournament Sidorenko property. One of the primary objectives of this section is to prove~\cref{c:transminusedge}, understanding when a transitive tournament minus an edge has the tournament Sidorenko property. Along the way, we exhibit several additional novel constructions of tournament Sidorenko oriented graphs.

We begin by recalling some basic constructions and properties of tournament Sidorenko oriented graphs, starting with the more general directed Sidorenko property studied in~\cite{genFHMZ22}.

\begin{defn}
A digraph $D$ is said to have the \textit{directed Sidorenko property} if for every digraph $G$,
\begin{equation}\label{e:sid1}
t(D, G) \ge t(\vec K_2, G)^{e(D)}.
\end{equation}
\end{defn}

In Theorem 1.5 of~\cite{genFHMZ22}, it is observed that if a digraph has a homomorphism to an edge, then it has the directed Sidorenko property if and only if its underlying undirected graph $\overline{D}$ satisfies the asymmetric Sidorenko property. Since a digraph $D$ has the tournament Sidorenko property if it has the \textit{directed Sidorenko property}, this immediately implies that a large class of digraphs have the \textit{tournament Sidorenko property}.

\begin{thm}[\cite{genFHMZ22}]
If $\overline{D}$ is asymmetric Sidorenko and $D$ has a homomorphism to an edge, then $D$ has the tournament Sidorenko property.
\end{thm}

We make a few further simple notes about digraphs with the tournament Sidorenko property.

\begin{obs} 
If both $D_1$ and $D_2$ have the tournament Sidorenko property, then so does their disjoint union $D = D_1 \sqcup D_2$.
\end{obs}

As noted in the introduction, transitivity is a prerequisite for a digraph to have the tournament Sidorenko property, as there are no copies of any non-transitive digraph $D$ in an $n$-vertex transitive tournament.

\begin{obs}\label{t:musttrans}
Any tournament Sidorenko digraph $D$ must be transitive.
\end{obs}

We establish the directed Sidorenko property for several new families of digraphs, beginning by noting the Sidorenko analogue of~\cref{t:antiextend}.

\begin{prop}\label{p:sidextend}
Let $D_1, D_2$ be a pair of tournament Sidorenko digraphs such that there is an independent subset $I \subset V(D_1)$ with $|I| = v(D_2)$. If all vertices in $I$ each have the same in- and out-neighborhood, then the digraph $D$ obtained by embedding $D_2 \hookrightarrow I$ is also tournament Sidorenko.
\end{prop}
\begin{proof}
This follows by an identical argument to~\cref{t:antiextend}, reversing the relevant inequalities.
\end{proof}

The following corollary of~\cref{p:sidextend} shows that the \textit{join} of two tournament Sidorenko digraphs (given by connecting every $v \in D_1, w \in D_2$ by an edge $(v, w)$) is tournament Sidorenko.
\begin{cor}
If $D_1, D_2$ are tournament Sidorenko, so is the \textit{join} of $D_1, D_2$, the digraph $D$ with $V(D) = V(D_1) \cup V(D_2), E(D) = E(D_1) \cup E(D_2) \cup \{(v_1, v_2) : v_1 \in D_1, v_2 \in D_2\}$
\end{cor}
\begin{proof}
We apply~\cref{p:sidextend} twice. Let $D_0 = K_{v(D_1), v(D_2)}^{\rightarrow}$ be the orientation of the complete bipartite graph with a homomorphism to an edge. Then $D_0$ is tournament Sidorenko and thus by~\cref{p:sidextend}, so is $D' = (V(D_0), E(D_0) \cup E(D_1))$ for any embedding of $D_1$ into the part of $D_0$ with $v(D_1)$ vertices. Applying~\cref{p:sidextend} again, this time embedding $D_2$ into the size $v(D_2)$ independent set of $D'$ (from the original bipartition) then gives the desired result.
\end{proof}

We give another method of leveraging known tournament Sidorenko digraphs to construct larger tournament Sidorenko digraphs, a special case of which is~\cref{c:growspecial}.

\begin{defn}
Given two digraphs $D_1, D_2$ and a distinguished vertex $v^* \in D_1$, the digraph obtained by \textit{substituting $D_2$ for $v^*$} is $D' = (V, E)$ where $V = V(D_1) \backslash \{v^*\} \sqcup V(D_2)$ and 
$$(u, v) \in E(D) \iff \begin{cases}
(u, v) \in E(D_1) &  u,v \in V(D_1) \backslash \{v^*\} \\
(u, v) \in E(D_2) & u,v \in V(D_2) \\
(u, v^*) \in E(D_1) & u \in V(D_1), v \in V(D_2) \\
(v^*, v) \in E(D_1) & u \in V(D_2), v \in V(D_1)
\end{cases}
$$
\end{defn}

The following somewhat complicated looking result gives a class of larger tournament Sidorenko digraphs constructed from smaller such digraphs (a special case of the below result is~\cref{c:growspecial}).
\begin{thm}\label{p:growtsid}
For two digraphs $D_1$ and $D_2$ with at least $2$ vertices which have the tournament Sidorenko property, and a vertex $v^*\in D_1$ such that $D_1[V(D_1)\backslash\{v^*\}]$ is tournament anti-Sidorenko, the digraph formed by substituting $D_2$ for $v^*$ is tournament Sidorenko. 
\end{thm}
\begin{proof}
We follow a similar approach as in the proof of~\cref{t:antiextend}.
Suppose that $D_1, D_2$ are digraphs with the tournament Sidorenko property. Let $D$ be the digraph formed by substituting $D_2$ for $v^* \in D_1$ chosen such that $D_1' := D_1[V(D_1) \backslash \{v^*\}]$ is tournament anti-Sidorenko. Let $T = (V, E)$ be a tournament on $n$ vertices.
Since $D_1'$ is tournament anti-Sidorenko,
$$N_L(D_1', T) \le 2^{-e(D_1')}n^{v(D_1')} = 2^{-e(D_1) + \deg(v^*)}n^{v(D_1)-1}.$$
Let $S$ be a $v(D_1')$-tuple of vertices in $V$ such that $T[S] \sim D_1'$. Given such an $S$, we let $U_S$ denote the set of vertices $v \in V$ such that $T[(S, v)] \sim D_1$ i.e. that together with $S$, form a labeled copy of $D_1$. Note that $N_L(D_1, T) = \sum_{S \mid T[S] \sim D_1'} |U_S|.$
With this characterization, we then have the following:
\begin{align*}
N_L(D, T) &= \sum_{S \mid T[S] \sim D_1'} N_L(D_2, T[U_S]) \\
&\ge (1 - o(1))2^{-e(D_2)} \sum_{S \mid T[S] \sim D_1'} |U_S|^{v(D_2)} \tag{$D_2$ t.-Sid.}\\ 
&\ge  (1 - o(1))2^{-e(D_2)} \cdot N_L(D_1', T) \left(\frac{  \sum_{S \mid T[S] \sim D_1'} |U(S)| }{N_L(D_1', T)} \right)^{v(D_2)} \tag{convexity} \\
&=  (1 - o(1))2^{-e(D_2)} \cdot N_L(D_1', T)^{1 - v(D_2)} N_L(D_1, T)^{v(D_2)} \\
&\ge  (1 - o(1))2^{-e(D_2)} \cdot \left(2^{-e(D_1')}n^{v(D_1')} \right)^{1-v(D_2)} \left(2^{-e(D_1)}n^{v(D_1)}\right)^{v(D_2)} \tag{$D_1$ t.-Sid., $D_1'$ t.-anti-Sid.} \\
&=  (1 - o(1))2^{-e(D)}n^{v(D)},
\end{align*}
recalling that $v(D) = v(D_1) + v(D_2) - 1$ and $e(D) = e(D_2)+e(D_1)+\deg(v^*)(v(D_2)-1)$.
\qedhere
\end{proof}

To prove~\cref{c:transminusedge}, we will leverage~\cref{p:sidextend} and~\cref{c:growspecial} in conjunction with establishing the Sidorenko property for a new family of digraphs $D_k$.

\begin{defn}
Define digraph $D_k$ to be the digraph on $k+2$ vertices with $V(D_k) = \{a, b, v_1, \ldots, v_k\}$ and $E(D_k) = \bigcup_{i = 1}^k \{(a, v_i), (v_i, b)\}$.
\end{defn}

Note that $D_1$ is a directed $2$-edge path and is tournament anti-Sidorenko. We show that for $k \ge 2$, $D_k$ is tournament Sidorenko.

An \textit{impartial digraph} $I$ with $k$ vertices has the property that the number of copies of $I$ in any $n$-vertex tournament $T$ is only a function of $n$ and $k$ (i.e. does not depend on the orientation of the particular tournament in question). The set of impartial digraphs was completely characterized in~\cite{ZHAO19}, but for our bounds here, it will suffice to note that there exists (up to isomorphism and/or reversing all edges), a single impartial digraph with $4$ vertices, the tree $\vec I_4$ with vertex set $\{a,b,c,d\}$ and edge set $\{ (a, b), (c, b), (d, c)\}$. In the proof of the following lemma, we only need that $N_L(\vec I_4, T) = \frac{1}{8} n^4 + o(n^4)$ for any $n$-vertex tournament $T$. Equipped with this fact, we first observe $D_2$ is tournament Sidorenko.

\begin{lemma}\label{l:d2}
$D_2$ is tournament Sidorenko. 
\end{lemma}
\begin{proof}
Let $T = (V, E)$ be a tournament on $n$ vertices.
Let $\vec P_3$ be the directed path with $3$ vertices and let $\vec I_4$ be an impartial, oriented tree with $4$ vertices. We use $\ind_{xy}$ to denote the indicator that $(x, y) \in E$.
We observe that 
\begin{align*}
0 &\le \sum_{x, z \in V} \left( \sum_{y \in V} \left(\ind_{xy} - \frac12\right) \left(\ind_{yz} - \frac12\right) \right)^2 \\
&= \sum_{w,x,y,z \in V} \left(\ind_{xy} - \frac12\right) \left(\ind_{yz} - \frac12\right) \left(\ind_{xw}- \frac12\right) \left(\ind_{wz} - \frac12\right) \\
&= \sum_{w,x,y,z \in V} \ind_{xy}\ind_{yz}\ind_{xw}\ind_{wz} -\frac12 \cdot 2 \cdot \sum_{w,x,y, z \in V} (\ind_{xw} \ind_{xy} \ind_{yz} + \ind_{xy} \ind_{yz} \ind_{wz}) \\
&\qquad + \frac14 n \left(N_L(\vec K_{1,2}, T) + N_L(\vec K_{2,1}, T) + 2 N_L(\vec P_3, T)\right) + \frac14 \cdot 2 \cdot  N_L(\vec K_2, T)^2 \\
&\qquad - \frac18 \cdot 4 \cdot n^2 \cdot N_L(\vec K_2, T) + \frac{1}{16} n^4 \\
&= N_L(D_2, T) - \frac12 \cdot 4 \cdot N_L(\vec I_4, T) + \frac14 \left( n \cdot N_L(\overline{P_3},K_n) + 2 \binom{n}{2}^2 \right) - \frac18 \cdot 4 n^2 \binom{n}{2} + \frac{1}{16} n^4\\
&=  N_L(D_2, T) - \frac{1}{16} n^4 + o(n^4).
\end{align*}
The final equation is obtained by substituting $N_L(\vec I_4, T) = \frac{1}{8} n^4 + o(n^4)$ and $N_L(\overline{P_3}, K_n) = n^3 + o(n^3).$ 
This implies the desired result.
\end{proof}

\begin{lemma}\label{l:dsid}
$D_k$ is tournament Sidorenko for $k > 2$.
\end{lemma}
\begin{proof}
Let $T = (V, E)$ be a tournament on $n$ vertices.
By applying H\"older's inequality and~\cref{l:d2} we have that for $k > 2$,
\begin{align*}
\left(\frac{1}{2^4} + o(1)\right) n^4 &\le N_L(D_2, T)  = \sum_{x, z \in V} \left(\sum_{y \in V} \ind_{xy} \ind_{yz} \right)^2 \\
&\le \left( \sum_{x, z \in V} \left( \sum_{y \in V} \ind_{xy} \ind_{yz}  \right)^k \right)^{\frac{2}{k}} \left( n^2 \right)^{\frac{(k-2)}{k}} = N_L(D_k, T)^{\frac{2}{k}} n^{\frac{2(k-2)}{k}}.\end{align*}
Rearranging gives the desired inequality.
\end{proof}

We are now ready to prove~\cref{c:transminusedge}. The condition in the statement is necessary in at least the case $k = 3$, as seen by considering the tournament anti-Sidorenko path on 3 vertices $\vec P_3$ that is one edge away from the $3$-vertex transitive tournament.

\begin{proof}[Proof of~\cref{c:transminusedge}]
Let $T_k$ be the transitive tournament with vertex set $[k]$ such that $E(T_k) = \{(i, j) : i < j\}$ and choose some $e = (i, j)$ with $j > i, j - i \neq 2$. By~\cref{l:d2} and~\cref{l:dsid} and noting that $D_0$ is a single directed edge, $D_{j - i - 1}$ is tournament Sidorenko since $j - i - 1 \neq 1$. Further by applying~\cref{p:sidextend}, the digraph obtained by embedding a copy of the transitive tournament $T_{j - i - 1}$ into the size $j - i - 1$ independent set of $D_{j - i -1}$ is tournament Sidorenko since transitive tournaments are tournament Sidorenko per~\cite{COR17}. Call the resulting tournament Sidorenko digraph on $j - i - 1$ vertices $A$. We then iteratively apply~\cref{c:growspecial}, adding $i-1$ vertices one by one to $A$ that are each complete to the current digraph (including previously added vertices) with in-degree $0$ and then one-by-one adding $n - j$ vertices to $A$ that are each complete to the current oriented digraph, each with out-degree $0$. The resulting digraph is exactly $T_k \backslash \{e\}$ and has the tournament Sidorenko property by construction.
\end{proof}

\section{Orienting trees}\label{s:trees}
\subsection{Stars} We prove \cref{t:star} by first verifying that the given orientations of an undirected star in the theorem statement do indeed have the tournament Sidorenko or anti-Sidorenko property (respectively). We subsequently demonstrate that all other orientations of an undirected star have neither the tournament Sidorenko property nor the tournament anti-Sidorenko property.

We first show that a balanced even star satisfies the strongly tournament anti-Sidorenko property.
\begin{prop}\label{p:even-star}
    Let $G$ be a directed star on $2k+1$ vertices with a central vertex $v$ whose in-degree and out-degree are both $k$. Then $(\{v\}, G)$ is strongly tournament anti-Sidorenko.
\end{prop}
\begin{proof}
Let $T$ be any tournament on $n$ vertex and $u$ be any vertex in $T$. Let the in-degree and the out-degree of $u$ in $T$ be $d_1$ and $d_2$ respectively. Then the number of labeled copies of $G$ in $T$ with $v$ mapped to $u$ is at most
\[d_1^kd_2^k \le \left(\frac{d_1+d_2}{2}\right)^{2k} = \left(\frac{n-1}{2}\right)^{2k} < 2^{-e(G)}n^{v(G)-1}.\qedhere\]
\end{proof}
\begin{proof}[Proof of~\cref{t:star}]
Let $S = (\{v^*, v_1, \ldots, v_s\}, E)$ be some orientation of a star, with $\overline{S} = K_{1, s}$ such that $v^*$ has total degree $s$ with $d^+ = d^+(v^*), d^- = d^-(v^*)$. We first note that if $S$ has a homomorphism to a directed edge, $S$ has the directed Sidorenko property per~\cite{genFHMZ22} and thus the tournament Sidorenko property.

We next consider the case where $|d^+ - d^-| \le 1$ and verify that in this case, $S$ has the tournament anti-Sidorenko property.
In fact, this is a consequence of \cref{l:glue-anti-sid}, where $(I_1, D_1)$ is $(\{v\}, G)$ as in \cref{p:even-star}, and $I_2=\emptyset$ and $D_2$ is either a single vertex (when $s$ is even) or a digraph with two vertices and an edge (when $s$ is odd). In the latter case, we can choose which vertex $v$ is identified with depending on whether $d^+-d^- = 1$ or $-1$.

We now show that in all other cases, $S$ does not have either the tournament Sidorenko nor the tournament anti-Sidorenko property. Assuming that we do not fall into the above two cases, we know that $d^+, d^- > 0$ with $|d^+ - d^-| \ge 2$. Without loss of generality, let $d^+ < d^-$. For $c \in [0, 1]$, we consider the family of random tournaments $T_n^{(c)} = ([n], E)$, where we have $(i, j) \in E$ if $i \le cn$ and $j > cn$, else we uniformly at random sample one of $(i, j), (j, i) \in E$. 
We count copies of $S$ in $T_n^{(c)}$ by considering two cases, whether the copy of its center $v^*$ is among the first $cn$ vertices or not. Viewing $S$ as a labeled digraph, this gives:
\begin{align*}
\EE\left[t_S\left(T_n^{(c)}\right)\right] &= c\left(1 - \frac{c}{2}\right)^{d^+} (c/2)^{d^-} + (1-c) \left(c + \frac{1-c}{2}\right)^{d^-}\left(\frac{1-c}{2}\right)^{d^+} \\
&= 2^{-s} \left(c^{1 + d^-} (2 - c)^{d^+} + (1-c)^{1 + d^+}(1 + c)^{d^-}\right) \\
&=: 2^{-s} f(c).
\end{align*}
Note that $f$ is a polynomial in $c$ with $f(0) = f(1) = 1$. Since $d^+, d^- > 0$, we have that
$f'(0) = d^- -d^+ - 1 > 0,$ since $d^+ \le d^- - 2$. Thus, by choosing a deterministic family of tournaments $T^+_n$ that achieves at least the above expectation and a suitable $c > 0$, we see that $S$ is not tournament anti-Sidorenko. Similarly, note that $f'(1) = d^- + 1 - d^+ \ge 1$, and thus the derivative at $1$ is also positive. Choosing a deterministic family of tournaments $T^-_n$ that achieves at most the above expectation for a suitable choice of $c < 1$ where $f(c) < f(1) = 1$, we see that $S$ as also not tournament Sidorenko. Thus, $S$ neither has the tournament Sidorenko nor the anti-Sidorenko properties.
\end{proof}

\subsection{General trees} 
Recall that~\cref{c:treeanti} posits that every undirected tree has an orientation which is tournament anti-Sidorenko. Below, we prove~\cref{p:partialtree}, which proves the conjecture for a specific family of trees. We then study the $1$-subdivisions of an even star $K_{1, 2k}$, observing that a natural choice of orientation of this graph is \textit{not} tournament anti-Sidorenko.

We prove the following strengthening of \cref{p:partialtree} for trees with one vertex of even degree.
\begin{prop}
    Let $R$ be an undirected tree with exactly one vertex $\rho$ of even degree. Then there exists an orientation $D$ of $R$ such that $(\{\rho\}, D)$ is strongly tournament anti-Sidorenko.
\end{prop}
\begin{proof}
Let $R$ be an undirected tree with vertex $\rho$ the sole vertex of even degree in $R$. We view $R$ as a \textit{rooted tree} with root $\rho$.
Then we know that every vertex in this rooted tree has an even number of children.
We show the existence of such an orientation $D$ by an induction on the number of vertices in $R$.

If $|V(R)| = 1$, then $R$ is a single vertex. The statement trivially holds.

Suppose $|V(R)| > 1$. The induction hypothesis is that every undirected tree $R'$ with $|V(R')| < |V(R)|$ which is rooted at a vertex $\rho'$ that is the unique vertex in $R'$ of even degree has an orientation $D'$  such that $(\{\rho'\}, D')$ is strongly tournament anti-Sidorenko.

Suppose that $R$ has depth $\ell$ when viewed as a tree rooted at $\rho$. Since $|V(R)| > 1$, we know that $\ell \ge 1$ (say the root $\rho$ has depth $0$). Let $v$ be a vertex of $(R, \rho)$ at depth $\ell - 1$ that has at least one child. Since every vertex has an even number of children, $v$ has at least two children, both of which are at depth $\ell$ and therefore are leaves. Let $x$ and $y$ be two children of $v$. Let $R'$ be the tree derived from removing $x$ and $y$ from $R$.

Note that $V(R') = V(R)\setminus \{x, y\}$ and $E(R') = E(R) \setminus \{\{v, x\}, \{v, y\}\}$. Hence $R'$ is also a tree with $\rho$ being the only vertex of even degree and has strictly fewer vertices than $R$. By the induction hypothesis there exists an orientation $D'$ of $R'$ such that $(\{\rho\}, D')$ is strongly tournament anti-Sidorenko. Now let $D_1$ be the digraph of three vertices $\{v, x, y\}$ and two edges $\{(x, v), (v, y)\}$. By \cref{p:even-star} with $k = 1$, we know that $(\{v\}, D_1)$ is strongly tournament anti-Sidorenko. By \cref{l:glue-anti-sid} with $(I_1, D_1) = (\{v\}, D_1)$ and $(I_2, D_2) = (\{\rho\}, D')$, and $v$ in $D_1$ is identified $v$ in $D_2$, we see that the derived digraph $D$ is an orientation of $R$, and $(\rho, D)$ is strongly tournament anti-Sidorenko.
This completes the proof by induction.
\end{proof}
\cref{p:partialtree} follows for trees with one vertex of even degree because $(\{\rho\}, D)$ being strongly tournament anti-Sidorenko implies that $(\emptyset, D)$ is strongly tournament anti-Sidorenko, or, equivalently, $D$ is tournament anti-Sidorenko. Also note that a single edge is tournament anti-Sidorenko, with \cref{l:glue-anti-sid} we may also conclude that if every vertex in a tree is of odd degree, then there is a tournament anti-Sidorenko orientation of this tree. This completes the proof of~\cref{p:partialtree}.
\begin{rem}

The above leaves open~\cref{c:treeanti} for most trees. A very simple example of a tree not covered by~\cref{p:partialtree} is a $1$-subdivision of a $K_{1, 2k}$, which has $2k + 1$ vertices of even degree. Naively, noting that directed paths are tournament anti-Sidorenko~\cite{SSZ20}, we might hope that by orienting this tree as a union of $k$ directed paths of length $4$, we achieve a tournament anti-Sidorenko digraph. However, for $k = 2$ this intuition already breaks down. Even though the number of copies of such a digraph in a \textit{transitive tournament} is under the random bound, one can find a family of tournaments where the number of copies of this digraph exceeds the random bound.

There is nonetheless a simple orientation of a $1$-subdivision of $K_{1, 2k}$ for even $k$ that is tournament anti-Sidorenko; namely the orientation where all $4$ possible orientations of $2$-edge paths occur equally often amongst the $2$-edge paths that have the central node as one endpoint.
\end{rem}

\section{Tournament forcing}\label{s:forcing}
We show that tournament forcing digraphs must be systematically under or over-represented.

\begin{proof}[Proof of \cref{t:overunder}]
 Suppose, by way of contradiction, that $D$ has the tournament forcing property but does not have the tournament Sidorenko property or the tournament anti-Sidorenko property. This implies there is an $\epsilon>0$ and two families of tournaments, $\{T_{n}'\}_{n=1}^{\infty}$ and $\{T_{n}''\}_{n=1}^{\infty}$, where we assume without loss of generality that $T_n', T_n''$ have $n$ vertices, such that for all $n$ sufficiently large, 
\begin{align*}
    t_{D}(T'_{n})\leq (1-\epsilon){2^{-e(D)}}, & & t_{D}(T''_{n})\geq(1+\epsilon){2^{-e(D)}}.
\end{align*}
It follows from quasirandomness for tournaments (the analogous result for undirected graphs is in~\cite{FAN89}) that for every $\delta>0$ which is sufficiently small in $\epsilon$ and $D$, and for sufficiently large $n$, tournament $T_n'$ has vertex subsets $X_n, Y_n$ with $|X_n| = |Y_n| = \lfloor \delta n \rfloor$ and $e(X_n,Y_n)<\frac{1-\delta}{2}|X_n||Y_n|$. Let $S_n=X_n \cup Y_n$, so $|S_n|=2\lfloor \delta n\rfloor$. Note that for any vertex subset $S$ of a tournament $T$ on $n$ vertices, the number of homomorphisms of a digraph $D$ to $T$ that maps to at least one vertex from $S$ is at most $v(D)|S| n^{v(D)-1}$. 

We suppose that the vertices of each $T_n', T_n''$ are each labeled by $[n]$ (in some, not necessarily consistent order). Let $m = \binom{n}{2} - \binom{|S_n|}{2}$ be the number of (unordered) pairs of distinct vertices that do not contain any vertices from $S_n$. Arbitrarily label these unordered vertex pairs as $p_1, \ldots, p_m$. 
We then iteratively define a sequence of tournaments starting with $T_n^{(0)} := T_n'$. For each $i \in [m]$, we initialize $T_n^{(i)} = T_n^{(i-1)}$. We consider the pair of vertices $p_i$, and let the edge supported on $p_i$ by $T_n^{(i)}$ be oriented to match the orientation in $T_n''$ (i.e. we possibly flip the orientation of the edge given by the vertices in $p_i$). By construction, at each step, 
$$|h_{D}(T_{n}^{(i+1)})-h_{D}(T_{n}^{(i)})|\leq n^{v(D)-2}.$$  
Further, note that all edges in $T_n^{(m)}$ have the same orientation as edges in $T_n''$ with the possible exception of a subset of the edges that are incident to at least one vertex from $S_n$. Since there are at most $\delta v(D) n^{v(D)}$ homomorphisms from $V(D) \rightarrow V(T_n'')$ that have as image at least one vertex from $S_n$, we see that 
$$h_D(T_n^{(m)}) \ge h_D(T_n'') - \delta v(D) n^{v(D)} \ge (1 + \eps/2) 2^{-e(D)} n^{v(D)}.$$ 

Since $h_{D}(T_{n}^{(0)})=h_{D}(T_{n}')\leq (1-\epsilon)2^{-e(D)}n^{v(D)}$, we know that there exists $M\leq m$ such that for $T_{n}:=T_{n}^{(M)}$, $h_{D}(T_{n})=(1+o(1))2^{-e(D)}n^{v(D)}$, and thus $t_{D}(T_{n})=(1+o(1))2^{-e(D)}.$ Since $D$ has the tournament forcing property, this implies that $\{T_{n}\}_{n=1}^{\infty}$ has quasirandom direction. However, this is impossible as $T_n[S_n] \cong T_n'[S_n]$, and the sets $S_n$ certify the lack of quasirandomness of the sequence $\{T_n'\}$ and therefore, certify that the sequence $\{T_n\}$ is also not quasirandom. This produces a contradiction, and we therefore conclude that $D$ is not forcing.
\end{proof} 

The above result has several immediate consequences; we observe a few below.
\begin{cor}
If non-transitive digraph $D$ is tournament forcing, then it is tournament anti-Sidorenko.
\end{cor}
\begin{proof}
If $D$ is non-transitive, then there are no copies of $D$ in any transitive tournament. Therefore, $D$ cannot have the tournament Sidorenko property, and thus by~\cref{t:overunder}, $D$ must be tournament anti-Sidorenko.
\end{proof}

This immediately implies~\cref{t:forcingntourn}. We also observe that one of our earlier constructions of a tournament anti-Sidorenko digraph has the tournament forcing property:

\begin{cor} 
For $k \ge 2$, the digraph $D$ in~\cref{p:anti2} is tournament forcing.
\end{cor}
We recall the digraph from \cref{p:anti2}: let $D=(A\sqcup B,F)$ be an orientation of the complete bipartite graph with parts of size $|A|=k,|B|=2^{k}$ such that for every subset $A'\subset A$ there is a unique vertex $v$ with $N^{+}(v)=A'$. 
\begin{proof}
Consider a sequence of $n$-vertex tournaments $\{T_n\}$ and let $T = T_n$.
The proof that $D$ was tournament anti-Sidorenko in~\cref{p:anti2} included one application of AM-GM. From this, we note that $N_L(D, T) = (1 + o(1)) 2^{-e(D)} n^{v(D)}$ if and only if equality holds in the AM-GM  up to a multiplicative $1 + o(1)$ factor; this happens exactly when all but a $o(1)$ fraction of the size $k$ subsets $S \subset V$ satisfy the following condition: for all $S' \subset S$,
$$\left| \{v \in V \backslash S \mid N^+(v) \cap S = S' \}\right| = (1 + o(1))(n-k)2^{-k}.$$
However, the above condition implies that $T$ has $\eps$-quasirandom direction (as defined in~\cite{GRI13}) for some $\eps = o(1)$ and thus by~\cite{FAN91,GRI13}, we find that $\{T_n\}$ is a quasirandom sequence of tournaments. 
\end{proof}

\section{Concluding remarks}\label{s:conc}
Counts of oriented subgraphs in tournaments exhibit a wide range of behaviors, not all of which have natural analogues in undirected graphs. In addition to the tournament anti-Sidorenko property we study here, there are a remarkable family of oriented forests, called \textit{impartial digraphs} (briefly discussed in proving~\cref{l:d2}), that have the property that the number of copies of a given impartial digraph in an $n$-vertex tournament is only a function of $n$ (i.e. independent of the orientation of the tournament). A single vertex is an impartial digraph, taking the disjoint union of impartial digraphs is an impartial digraph, and taking the disjoint union of two copies of an impartial digraph $D$ and adding an edge between the two copies of a single vertex from $D$ is an impartial digraph. Zhao and Zhou in~\cite{ZHAO19} proved a conjecture of Fox, Huang, and Lee that this gives a complete characterization of impartial digraphs. 

It would be very interesting to characterize the family of all tournament anti-Sidorenko digraphs, and the family of tournament Sidorenko digraphs. It should be noted that, if a digraph is both tournament anti-Sidorenko and tournament Sidorenko, then it is an impartial digraph \cite[Proposition 3.2]{ZHAO19}. This means that we understand the intersection of the two families. Short of a complete characterization of either family mentioned above, the following pair of questions would be interesting to answer. First, we might hope to resolve~\cref{c:treeanti}.

\begin{qn}
Does every undirected tree have a tournament anti-Sidorenko orientation?
\end{qn}

Note that tournament anti-Sidorenko digraphs have to be quite sparse per~\cref{c:notanti}. However, tournament Sidorenko graphs can be very dense (e.g. transitive tournaments have the tournament Sidorenko property). This can be rephrased as saying that $K_n$ has a tournament Sidorenko orientation,~\cref{c:transminusedge} implies that $K_n$ minus an edge has a tournament Sidorenko orientation, and~\cite{genFHMZ22} implies that every undirected graph with the asymmetric Sidorenko property has a tournament Sidorenko orientation. These examples suggest that it might be interesting to study which undirected graphs have some tournament Sidorenko orientation, beginning with dense graphs.

\begin{qn}
Does every undirected graph with $k$ vertices and $\Omega(k^2)$ edges have a tournament Sidorenko orientation?
\end{qn}

Fully understanding the tournament forcing property may be even more challenging. It is not necessary for a digraph to be transitive to have the tournament forcing property. Some additional conditions are necessary, as the tournament forcing, non-transitive oriented tournament $T_5^8$ in Figure 2 of~\cite{COR15} illustrates. The following question remains open.

\begin{qn}
Which tournament anti-Sidorenko and Sidorenko digraphs have the tournament forcing property?
\end{qn}

\bibliographystyle{alpha}
\bibliography{refs.bib}

\clearpage 
\appendix

\end{document}